\documentclass[a4paper,12pt]{article}
\usepackage{times, url}
\usepackage{amsmath,amssymb,amsthm,amsfonts}

\newtheorem{theorem}{Theorem}[section]

\newtheorem{proposition}[theorem]{Proposition}

\theoremstyle{definition}
\newtheorem{definition}{Definition}[section]

\theoremstyle{remark}

\begin{document}
\setcounter{page}{1}

\begin{center}
{\LARGE \bf  On sequences of natural numbers having pairwise relatively prime terms.}
\vspace{8mm}

{\large \bf Gaitanas Konstantinos}
\vspace{3mm}

 Department of  Applied Mathematical and Physical Sciences\\
 National Technical University of Athens \\ 
Heroon Polytechneiou Str., Zografou Campus, 15780 Athens, Greece \\
e-mail: \url{raffako@hotmail.com}
\vspace{2mm}

\end{center}
\vspace{10mm}

\noindent
{\bf Abstract:}We prove some theorems which give sufficient conditions for the existence of prime numbers among the terms of a sequence which has pairwise relatively prime terms.\\

\vspace{10mm}

\section{Introduction} 

There are many known sequences of numbers having all terms pairwise relatively prime and their study has important applications for the theory of  numbers.Goldbach  in a letter to Leonard Euler in order to prove that there are infinitely many primes used such a sequence observing that all of its terms have distict prime divisors.Nowadays, we could say that some of the most studied examples of such sequences are the Fermat numbers $F_n=2^{2^n}+1$, numbers of the form $2^p-1$ with $p$ being prime, or Fibonacci numbers with a prime index.\\It has been conjectured that the above three sequences contain infinitely many primes but this still remains an open problem.\\
We will not study the behaviour of these sequences. Instead we will obtain more generall results which shall give new insight in this type of sequences and prove that it is possible to obtain prime values if certain conditions are satisfied.\\

\section{Preliminaries}
Throughout this paper, $\pi(n)$ stands for the prime counting function, $P$ denotes the set of prime numbers and  $P(s)=\displaystyle\sum\limits_{p} \frac{1}{{p}^s}$ denotes the prime zeta function.\\
Especially we are interested in $P(2)=\displaystyle\sum\limits_{p} \frac{1}{{p}^2}\approx0.4523$.\\
\begin{definition}We say that a sequence of natural numbers is a  PLP sequence if is strictly increasing and its terms are pairwise relatively prime and greater than $1$.
\end{definition}
  \section{Main results}
All the theorems we are going to see describe the behaviour of pairwise relatively prime numbers.We will prove that a sequence which is PLP and contains only composites  grows ``too fast'', a fact which plays a key role in order to obtain the following results.
\begin{proposition}
Let $a_n$ be a PLP sequence having all of its terms composite.\\
Then the series $\displaystyle\sum\limits_{n=1}^\infty \frac{1}{a_n}$ converges.
\begin{proof}
Let $p_n$ denote the least prime divisor of $a_n$. Since $a_n$ is composite we can write $a_n=p_n\cdot k$ with $k\geq p_n$.This means $a_n\geq {p_n}^2\Rightarrow\frac{1}{a_n}\leq\frac{1}{{p_n}^2}$. Since all terms are pairwise relatively prime all the prime divisors are distinct and thus the series is not greater than the sum of the reciprosals of the squares of all primes.\\
All these give  $\displaystyle\sum\limits_{n=1}^\infty \frac{1}{a_n}\leq P(2)$ so it is clear that the series converges to a sum not greater than $0.4523$.
\end{proof}
\end{proposition}
The above series is equal to $P(2)$ if we let for every $n\geq 1$ , $a_n={p_n}^2$.\\
 The following theorem is an immediate consequence of proposition 3.1
\begin{theorem}
Let $a_n$ be a PLP sequence.
Then if $\displaystyle\sum\limits_{n=1}^\infty \frac{1}{a_n}=\infty $ the sequence contains infinitely many prime numbers.
\begin{proof}
Suppose that from some point on, $a_n$ produces only composites.\\
We shall reach the absurd conclusion that the series $\displaystyle\sum\limits_{n=1}^\infty \frac{1}{a_n}$ is in fact bounded.\\
Let $k\in \mathbb{N}$ and suppose that for every $n\geq k$,  $a_n$ is composite.From theorem 3.1 we can see that\\
$\displaystyle\sum\limits_{n=1}^\infty \frac{1}{a_n}=\displaystyle\sum\limits_{n=1}^{k-1} \frac{1}{a_n}+\displaystyle\sum\limits_{n=k}^\infty \frac{1}{a_n}\leq \displaystyle\sum\limits_{n=1}^{k-1} \frac{1}{a_n}+P(2)<\infty$,  since the sum $\displaystyle\sum\limits_{n=1}^{k-1} \frac{1}{a_n}$ is bounded for fixed $k$.\\
\\
This contradicts the hypothesis that the above series diverges and thus theorem 3.2 holds true.
\end{proof}
\end{theorem}
No non trivial sequence satisfying the hypothesis of theorem 3.1 is known.Of course the sequence of prime numbers itself is such a sequence and also every such sequence which has asymptotic density among the primes.It would be desirable to construct such a sequence, but this construction is well beyond  the reach of known methods.\\
We can observe that actually ``almost all" of the terms of $a_n$ in theorem 3.1 are primes.If there are composites, proposition 3.1 implies that the sum of their reciprocals converges and thus the sum $\sum\limits_{a_n\in P} \frac{1}{a_n}$ diverges.\\
The following proposition is actually a generalization of the previous theorem.

\begin{proposition}Let $a_n$ be a sequence which has the PLT property.\\
Then if $\displaystyle\sum\limits_{n=1}^\infty \frac{1}{a_n}>P(s)$ there is a $k$ for which $a_k$ has at most $s-1$ prime factors.\\
Especially, if $\displaystyle\sum\limits_{n=1}^\infty \frac{1}{a_n}>P(2)$ there is at least a natural number $k$, for which $a_k$ is prime.
\begin{proof}
The proof follows the lines of theorem 3.1.\\
If every term had at least $s$ prime divisors then $a_n\geq {p_n}^s$ which gives $\displaystyle\sum\limits_{n=1}^\infty \frac{1}{a_n}\leq P(s)$ which leads to an immediate contradiction.
\end{proof}
\end{proposition}
It is reasonable to ask how dense must a PLP sequence be in order to contain prime values.We will show that such a sequence which has $\pi(\sqrt{n})+1$ terms contains at least one prime number and we will make use of a convenient approximation proved by Rosser and Schoenfeld in order to obtain a more accurate result.

\begin{proposition}Every PLP sequence which has $\frac{2\sqrt{n}}{\ln n}\cdot (1+\frac{3}{\ln n})+1$ terms not exceeding $n$ contains at least one prime number.\\
\begin{proof}
Suppose we choose $k$ pairwise relatively prime numbers $a_1,\ldots , a_k$ not greater than $n$.Let $p_i$ denote the least prime divisor of  $a_i$.If all the terms of the sequence are composite, we may proceed in a similar way to the proofs of the previous theorems:\\
For every $i$ with $1\leq i\leq k$,  \, $ {p_i}^2\leq a_i\leq n$ holds.This means all the prime divisors of the $a_i$'s are not greater than $\sqrt{n}$ which gives the bound $k\leq \pi(\sqrt{n})$.\\
It suffices to prove that the length of the sequence $k$ is in fact greater than $\pi(\sqrt{n})$ and conclude that all the terms cannot be composite.\\
It is proved \cite{1} that for every $n>1$, $\pi(n)<\frac{n}{\ln n}\cdot (1+\frac{3}{2\ln n})$ holds.Substituting $n$ by $\sqrt{n}$ yields  $\pi(\sqrt{n})<\frac{2\sqrt{n}}{\ln n}\cdot (1+\frac{3}{\ln n})$.\\
But from our assumption $k\geq\frac{2\sqrt{n}}{\ln n}\cdot (1+\frac{3}{\ln n})+1$ and we can see that all these together yield \\
$\frac{2\sqrt{n}}{\ln n}\cdot (1+\frac{3}{\ln n})+1\leq k\leq \pi(\sqrt{n})<\frac{2\sqrt{n}}{\ln n}\cdot (1+\frac{3}{\ln n})$ which is not possible.\\
Therefore at least one term of the sequence is prime and the proof is complete.
\end{proof}
\end{proposition}

\makeatletter
\renewcommand{\@biblabel}[1]{[#1]\hfill}
\makeatother

\end{document}